\documentclass[12pt, a4paper]{amsart}
\usepackage[utf8]{inputenc}
\usepackage{amsxtra,amssymb,amsthm,amsmath,amscd,mathrsfs, epsfig, eufrak}
\usepackage{amscd, amsmath, mathrsfs, amssymb, amsthm, amsxtra, bbding, epsfig, eucal, eufrak, graphicx, latexsym, mathrsfs, mathbbol, bbold}
\usepackage[all]{xy}

\usepackage{tikz}

\usepackage[normalem]{ulem}
\usepackage{soul}
\usepackage{color}

\setstcolor{red}

\makeatletter
\@namedef{subjclassname@2010}{%
\textup{2010} Mathematics Subject Classification}
\makeatother

\def\le{\leqslant}

\theoremstyle{plain}

\theoremstyle{remark}

\numberwithin{equation}{section}

\addtocounter{footnote}{1}

\begin{document}


\title[\tiny Sums of eight fourth power of primes]
{Sums of eight fourth power of primes}
\author[\tiny Yang Qu \& Rong Ma]{Yang Qu \& Rong Ma}

\address{%
Rong Ma
\\
School of Mathematics and Statistics
\\
Northwestern Polytechnical University
\\
Xi'an
\\
Shaanxi 710072
\\
China}
\email{marong0109@163.com}

\address{%
Yang Qu
\\
School of Mathematics and Statistics
\\
Northwestern Polytechnical University
\\
Xi'an
\\
Shaanxi 710072
\\
China}
\email{1993845352@qq.com}

\date{\today}

\begin{abstract}
For any sufficiently large $\ell$, suppose that $\ell$ can be  expressed as
$$ \ell=p_1^4+p_2^4+p_3^4+ \cdots +p_8^4,$$
where $p_1, p_2,p_3,\cdots, p_8$ are primes.For such $\ell$, in this paper we will use circle method and sieves to prove that the proportion of $\ell$ in positive integers is at least $\frac{1}{414.465}$ .
\end{abstract}

\keywords{Waring–Goldbach problem, Circle method, Sieve methods }

\maketitle

\section{Introduction}
For any positive integers $m$ and $k$, the Waring-Goldbach Problem is to discuss every positive number $\ell_{m,k}$ can be expressed as the sum of the $m$-$th$ powers of $k$ primes $p_1,p_2,\cdots,p_k$, i.e.
\begin{equation}
   \ell_{m,k}=p_1^m+p_2^m+\cdots p_{k-1}^m+p_k^m.
\end{equation}

When $m=1,k=2 $, this is the famous Goldbach's Conjecture, the difficult problem, which has not been successfully solved so far. But when $m=1$, $k=3$, the case of three prime numbers has been proven, which is the famous theorem of three prime numbers\textsuperscript{\cite{9}}.  Furthermore, for higher-order terms $m\geq1$, scholars hope to obtain the smallest $k$ that enables all positive integers to be expressed in the form of (1.1), and denote this smallest $k$ as $G(m)$. 

In 1770, Waring deduced from a finite set of evidence that every positive integer is the sum of four squares, nine cubes, nineteen fourth powers, and so on. In 1770, Lagrange proved the existence of $G(2)$, and in the following 139 years, proofs of existence were obtained for $m=3, 4, 5, 6, 7, 8,  10$. In 1909, Hilbert used induction to prove the existence of $G(m) $for $m$\textsuperscript{\cite{13}}.

In addition, for those $m,k$ that cannot express all positive integers in the form of (1.1), scholars also hope to study the density of positive integers that can be expressed in the form of (1.1) among all positive integers in this case. For example $m=3,k=4 $, some scholars have already obtained some results. Let $\mathcal{L}_{3,4}$ be the set of integers $\ell_{m,k}$ that can be written as the formula (1.1) when $m=3,k=4 $ . In 1949, Roth\textsuperscript{\cite{5}} showed that
$$\sum_{\ell_{3,4}\leqslant N\atop\ell_{3,4}\in\mathcal{L}_{3,4}}1\gg\frac N{\log^8N}.$$
Furthermore, from 2001 to 2003, Ren\textsuperscript{\cite{2}}\textsuperscript{\cite{3}}improved Roth’s theorem to the extent that
$$\sum\limits_{\substack{\ell_{3,4}\leqslant N\\\ell_{3,4}\in\mathcal{L}_{3,4}}}1\geqslant\beta_{3,4} N,$$
where $\beta_{3,4}=1/320$. In addition, Liu\textsuperscript{\cite{4}} improved this result and obtained $\beta_{3,4}=1/173.12$ in 2012.\\

The authors are interested in this problem especially for more cases $m$ and $k$. In this paper we consider the case $m=4,k=7 $ and try to say something about it. Let $\ell$ be a positive integer, suppose $\ell$ can be expressed as the sum of the fourth powers of eight prime numbers, i.e.
\begin{equation}
   \ell=p_1^4+p_2^4+p_3^4+\cdots +p_8^4,
\end{equation}
define $\mathcal{L}$ as the set of all $\ell$, and assume that
\begin{equation}
    \sum\limits_{\substack{\ell_\leqslant N\\\ell\in\mathcal{L}}}1\geqslant\beta N,
\end{equation}
Then we will consider the distribution of $\ell$. No one has yet conducted research on this problem, at least the authors have not see any relevant references. Therefore, we have the following theorem.\\

\noindent\textbf{Theorem.} \textit{$\beta=\frac{1}{414.465}$ is acceptable in (1.3).}\\

\noindent\textbf{Note.} Compared to the sums of four cubes of primes, the density of sums of eight fourth power of primes is smaller. For cases with a relatively small number, the results obtained by the method in this paper are not ideal, but if we expand the expression to the sums of more than eight prime numbers, the density may become even higher, which is a problem worthy of further research.
\section{Prerequisite Knowledge}
Firstly, we introduce some definitions to prepare for the proof of the theorem.\\

 Let $N$ be a large integer,  $\delta_0$ be small enough and
 $$U=\left(\frac{N}{64(1+\delta_0)}\right)^{1/4},\quad V=U^{7/8}.$$

Now, we define $A$ and $B $ as two sufficiently large integers, and $B $ is sufficiently large relative to $A $. Define $\mathfrak{M}(q,a)$  as the interval $[a/q-L^B/U^4,a/q+L^B/U^4]$, then write $\mathfrak{M}$ for the union of all $\mathfrak{M}(q,a)$ with $1\leqslant a\leqslant q\leqslant L^B$ and  $(a,q)=1.$ Obviously, $\mathfrak{M}(q,a)$ are disjoint. Define $\mathfrak{m}$ as the complement of $\mathfrak {M} $ in $[L ^ B/U ^ 4,1+L ^ B/U ^ 4] $.\\

For any real numbers $\lambda$ and $X$, define the following two integrals
$$\Phi(\lambda,X)=\int\limits_{X}^{2X}e(u^4\lambda)\:du,\quad\Psi(\lambda,X)=\int\limits_{X}^{2X}\frac{e(u^4\lambda)}{\log u}\:du.$$
Let
\begin{equation}
    J(n)=\int_{-U^{-4}L^B}^{U^{-4}L^B}\Phi(\lambda,U)\overline{\Psi(\lambda,U)}\Big|\Psi(\lambda,U)\Big|^6\Big|\Psi(\lambda,V)\Big|^8e(-n\lambda)\:d\lambda.
\end{equation}
Define
\begin{equation}
    \mathfrak{S}_{d}(n)=\sum_{q=1}^{\infty}T_{d}(n,q)\:,
\end{equation}
where
\begin{equation}
    T_d(n,q)=\sum_{\substack{a =1\\(a,q)=1}}^q\frac{S(q,ad^4)C^7(q,a)\overline{C^8(q,a)}}{q\varphi^{15}(q)}e\biggl(-\frac{an}{q}\biggr)
\end{equation}
and
\begin{equation}
    S(q,a)=\sum_{m=1}^qe\biggl(\frac{a m^4}{q}\biggr),\quad C(q,a)=\sum_{\substack{q=1\\(q,m)=1}}^qe\biggl(\frac{a m^4}{q}\biggr).
\end{equation}
\section{Some Lemmas}
To prove the validity of the theorem, we also introduce the following lemmas:\\

\noindent{\textbf{Lemma 1.}} \textit{The number S of solutions of 
\begin{equation}
    x_1^4+y_1^4+y_2^4=x_2^4+y_3^4+y_4^4
\end{equation}
with $U<x\leq2U,U^{7/8}<y\leq2U^{7/8}$ satisfies}
$$S\ll U^{25/8}.$$
\noindent{\textbf{Proof.}} For $x_1=x_2$, we can easily calculate that the number $S$ of solutions of (3.1) satisfies
$$S\ll U^{25/8}.$$
For $x_1\neq x_2$, base on symmetry, we assume $x_1<x_2$ and write $x_2=x_1+h$. Then (3.1) becomes
$$h(4x_1^3+6x_1^2h+4x_1h^2+h^3)=y_1^4+y_2^4-y_3^4+y_4^4.$$
Since $y_1^4+y_2^4\leq 32U^{7/2}$ and $x_1^3>U^3$ it follows that $h<8U^{1/2}$.

Let
$$F_{h}(\alpha)\:=\:\sum_{U<x\:\leqslant\:2U}e\bigl(\alpha h(4x_1^3+6x_1^2h+4x_1h^2+h^3)\bigr)\:,$$
$$G(\alpha)=\sum_{h<8U^{1/2}}F_{h}(\alpha)\:,$$
$$f(\alpha)\:=\:\sum_{U^{7/8}<y\:\leqslant\:2U^{7/8}}e(\alpha y^{4})\:.$$
Then 
$$S\leqslant\int_{0}^{1}G(\alpha)|f(\alpha)|^{4}d\alpha.$$
We hope to obtain an upper bound for $G(\alpha)$. Suppose that $|\alpha-a/q|<q^{-2}$ and $(a,q)=1$. By Cauchy' s inequality 
$$|G(\alpha)|^{2}\:\ll U^{1/2}\:\sum_{h\:<\:8U^{1/2}}|F_{h}(\alpha)|^{2}.$$
Moreover, write $y=x+h_1$ and we have
    \begin{align*}
        |F_{h}(\alpha)|^{2}\:&=\sum_{U<x\:\leqslant\:2U}\sum_{U<y\:\leqslant\:2U}e\bigl(\alpha h(4x_1^3+6x_1^2h+4x_1h^2-4y_1^3-6y_1^2h-4y_1h^2\bigr)\\&=\:\sum_{|h_1|<U}\sum_{\max(U,U-h_1)<x\leqslant\:\min(2U,2U-h_1)}e\bigl(\alpha h(12h_1x^2+12h_1(h+h_1)x+4h_1^3+6hh_1^2+4h^2h_1)\bigr).
    \end{align*}
Record the innermost sum formula as $T(h_1)$. Use Cauchy' s inequality again
$$|F_{h}(\alpha)|^{2}\:\leq U^{1/2}\biggl(\sum_{|h_1|<U}|T(h_1)|^2\biggl)^{1/2}.$$
Similarly, write $y'=x'+h_1$, then
\begin{align*}
    |T(h_1)|^2&=\sum_{|h_2|<U}\sum_{\max(U,U-h_2)<x'\leqslant\:\min(2U,2U-h_2)}e(\alpha hh_1h_2(24x+12(h+h_1+h_2)))\\&\ll U+\sum_{0<h_2<U}\min(U,||24\alpha hh_1h_2||^{-1}).
\end{align*}
Hence
\begin{align*}
    |F_{h}(\alpha)|^{2}&\ll U^{1/2}\biggl(\sum_{|h_1|<U}\bigg(U+\sum_{0<h_2<U}\min(U,||24\alpha hh_1h_2||^{-1})\bigg)\biggl)^{1/2}\\&\ll U^{3/2+\epsilon}+\sum_{0<u<24U^{2}h}\min(U,||\alpha u||^{-1}).
\end{align*}
By Lemma 2.2 of Vaughan\textsuperscript{\cite{7}}, we have
\begin{equation}
    |G(\alpha)|^{2}\ll U^{5/2+\epsilon}+U^{13/4+\epsilon}q^{-1/2}+U^{11/4+\epsilon}+U^{3/2+\epsilon}q^{1/2}
\end{equation}\\
Let $\mathfrak{M}'(q,a)$ denote the interval $[a/q-q^{-1}U^{-2},a/q+q^{-1}U^{-2}]$ and $\mathcal{U}=(P^{-2},1+P^{-2}]$. We may suppose that $U\geqslant 4$. Then the $\mathfrak{M}'(q,a)$ with $1\leqslant a\leqslant q\leqslant U$, $(a,q)=1$ are disjoint and contained in $\mathcal{U}$. Let $\mathfrak{M}$ be the union of the $\mathfrak{M}'(q,a)$ with $1\leqslant a\leqslant q\leqslant U$, $(a,q)=1$, and let $m'=\mathcal{U}\backslash \mathfrak{M}'$. Then
$$S\leqslant\:\int_{_{\mathfrak{M}'}}G(\alpha)|f(\alpha)|^{4}d\alpha\:+\:\int_{_{m'}}G(\alpha)|f(\alpha)|^{4}\:d\alpha.$$
By (3.2), when $q>U$, there is $\alpha \notin \mathfrak{M}'$ and 
$$G(\alpha)\ll P^{11/8+\epsilon}.$$
Hence, by Lemma 2.5 of Vaughan\textsuperscript{\cite{7}}, 
$$\int_{_{m'}}G(\alpha)|f(\alpha)|^{4}\:d\alpha\ll P^{25/8}.$$
For $\alpha \in \mathfrak{M}'$, there is $1\leqslant a\leqslant q\leqslant U$, then by Lemma 6.3 of Vaughan\textsuperscript{\cite{7}}
$$G(\alpha)\ll U^{13/8}q^{-1/4},$$
$$f(\alpha)\ll U^{7/8}q^{-1/4}(1+U^{7/2}|\alpha-a/q|)^{-1}.$$

Therefore
\begin{align*}
    \int_{_{\mathfrak{M}'}}G(\alpha)|f(\alpha)|^{4}d\alpha&\ll \sum_{q\leq U}\sum_{\stackrel{a=1}{(a,q)=1}}^{q}U^{13/8+7/2}q^{-1-1/4}\int_{\mathfrak{M}'(q,a)}\left(1+U^{7/2}\bigg|\alpha-\frac{a}{q}\bigg|\right)^{-4}d\alpha\\&\ll U^{19/8}
\end{align*}
This proves Lemma 1.\\

\noindent{\textbf{Lemma 2.}} \textit{Let $0\leqslant|n|\leqslant N$, for each $m$ with $U<m\leqslant2U$, denote by $R(m)$ the number of solutions of}
$$n=m^4+p_2^4+\cdots+p_8^4-p_9^4-\cdots-p_{16}^4$$
\textit{with}
$$p_2,p_3,p_4,p_9,p_{10}, p_{11},p_{12}\sim U,\quad p_5,p_6,p_7,p_{8},p_{13},p_{14},p_{15},p_{16}\sim V.$$
\textit{For $0 < \xi < 9/25$ and $D = N^\xi$, take $\mathfrak{S}_d(n)$ and $J(n)$ from the above (2.1), (2.2), define $E_d (n) $ as follows: }
$$\sum_{\substack{m \sim U \\ m\equiv0({\mathrm{mod}}~d)}} R(m) = \frac{\mathfrak{S}_d(n)}{d} J(n) + E_d(n).$$
\textit{Then we have:}\\
(i)\textit{$\mathfrak{S}_d(n)$ is absolutely convergent and satisfies} $\mathfrak{S}_{d}(n)\ll1.$\\
(ii)\textit{ $J(n)$ is positive and satisfies}
\begin{equation*}
    J(n)\leqslant KU^4V^8L^{-15}.
\end{equation*}
\textit{where $K=4888799.222$}

\noindent(iii)\textit{ For any complex numbers $\eta_{d}$ with $|\eta_{d}|\leqslant\tau(d)$, we hvae
\begin{equation*}
    \sum_{d\leqslant D}\eta_dE_d(n)\ll U^4V^8L^{-A}.
\end{equation*}\\}
\textbf{\noindent{Proof.}} Firstly, we prove (ii) of Lemma 2. From elementary estimation
$$\int_{X^4}^{16X^4}e(\lambda u)\:du\leqslant\min(X^4,|\lambda|^{-1})$$
and integration by parts, we have
$$\Phi(\lambda,X)=\frac{1}{4}\int_{X^4}^{16X^4}u^{-3/4}e(\lambda u)\,du\ll X^{-3}\min(X^4,|\lambda|^{-1})$$
and
$$\Psi(\lambda,X)\ll X^{-3}\log^{-1}X\min(X^4,|\lambda|^{-1}).$$

We hope to calculate (2.1) by integrating the entire real axis, but this will result in the following error
$$\begin{aligned}&\ll\int_{U^{-4}L^{B}}^{\infty}|\Phi(\lambda,U)||\Psi(\lambda,U)|^{7}|\Psi(\lambda,V)|^{8}\:d\lambda\\&\ll U^{-24}V^{-24}L^{-15}\int_{U^{-4}L^{B}}^{\infty}\min(U^{4},|\lambda|^{-1})^{6}\min(V^{4},|\lambda|^{-1})^{8}\:d\lambda\\&\ll U^4V^{8}L^{-5B}.\end{aligned}$$
By integral transformation, we have
$$\begin{aligned}
&\int_{-\infty}^{\infty}\Phi(\lambda,U)\overline{\Psi(\lambda,U)}\big|\Psi(\lambda,U)\big|^{6}\big|\Psi(\lambda,V)\big|^{8}e(-n\lambda)\:d\lambda\\
=&\frac{1}{4}\int_{\mathcal{D}}\frac{d\nu_{1}\cdots d\nu_{8}du_{1}\cdots du_{7}}{\nu_{1}^{3/4}\cdots\nu_{8}^{3/4}u_{1}^{3/4}\cdots u_{8}^{3/4}\log\nu_{1}\cdots\log\nu_{8}\log u_{1}\cdots\log u_{7}},
\end{aligned}$$
where
$$\mathcal{D}=\left\{(v_{1},\ldots,v_{8},u_{1},\ldots,u_{7})\colon V^{4}\leqslant v_{1},\ldots,v_{8}\leqslant 16V^{4},\ U^{4}\leqslant u_{1},\ldots,u_{7}\leqslant 16U^{4}\right\}$$
and
$u_{8}=n+v_{1}+v_{2}+v_{3}+v_{4}-v_{5}-v_{6}-v_{7}-v_{8}+u_{1}+u_{2}+u_{3}+u_{4}-u_{5}-u_{6}-u_{7}.$\\
Then we get
\begin{equation}
\begin{split}
J(n)=&\int_{-\infty}^{\infty}\Phi(\lambda,U)\overline{\Psi(\lambda,U)}|\Psi(\lambda,U)|^{6}|\Psi(\lambda,V)|^{8}e(-n\lambda)d\lambda\\
&+\:O(U^4V^{8}L^{-5B})\\
=&\frac{1}{4}\int_{V^{4}}^{16V^{4}}\:\frac{dv_{1}}{v_{1}^{3/4}\log v_{1}}\cdots\int_{V^{4}}^{16V^{4}}\:\frac{dv_{4}}{v_{8}^{3/4}\log v_{8}}\int_{U^{4}}^{16U^{4}}\:\frac{du_{1}}{u_{1}^{3/4}\log u_{1}}\\
&\times\int_{U^{4}}^{16U^{4}}\frac{du_{2}}{u_{2}^{3/4}\log u_{2}}\cdots\int_{\max(U^{4},x-16U^{4})}^{\min(16U^{4},x-U^{4})}\frac{du_{7}}{u_{7}^{3/4}(x-u_{7})^{3/4}\log u_{7}}\\
&+\:O(U^4V^{8}L^{-5B}).
\end{split}
\end{equation}
Let the last integral in the equation be denoted as $I $, where $x = n+v_{1}+v_{2}+v_{3}+v_{4}-v_{5}-v_{6}-v_{7}-v_{8}+u_{1}+u_{2}+u_{3}+u_{4}-u_{5}-u_{6}$, now we calculate the upper bound of $I $. Firstly, regarding $x\leqslant2U^{4}$ or $x\geqslant32U^{4}$, in both cases, the integral region does not exist, so we have $2U^{4}<x<32U^{4}$. In this case, we record as $u_7=xu $, then there is
\begin{align*}
    \begin{split}
        I&\leqslant(1+\varepsilon)x^{-1/2}L^{-1}\int_{\max(U^{4}/x,1-16U^{4}/x)}^{\min(16U^{4}/x,1-U^{4}/x)}u^{-3/4}(1-u)^{-3/4}\,du\\
&\leqslant\frac{1+\varepsilon}{\sqrt{2}}U^{-2}L^{-1}\int_{1/17}^{16/17}u^{-3/4}(1-u)^{-3/4}\,du\\&\leqslant\frac{1+\varepsilon}{\sqrt{2}}U^{-2}L^{-1}I ^ *,
    \end{split}
    \end{align*}
where we record as $I ^ *= \int_{1/17}^{16/17}u^{-3/4}(1-u)^{-3/4}\,du $, then substituting it into (3.3) yields the following
\begin{alignat*}{2}
    J(n)&\leq \frac{1+\varepsilon}{4\sqrt{2}}I^{*}U^{-2}L^{-1}\Bigg(\:\int_{U^{4}}^{16U^{4}}\:\frac{du}{u^{3/4}\log u}\Bigg)^{6}\Bigg(\:\int_{V^{4}}^{16V^{4}}\:\frac{dv}{v^{3/4}\log v}\Bigg)^{8}+O(U^4V^{8}L^{-5B})
    \\&<\left(\frac{32}{7}\right)^8\frac{4^5}{\sqrt{2}}(1+\varepsilon)I^*U^4V^8L^{-15}+O(U^4V^{8}L^{-5B})
    \\&<4888799.222 U^4V^{8}L^{-15},
\end{alignat*}
where we have used the estimate $I^{*}<7.73$. This proves (ii). 
Now we prove (iii), let
$$v_d(n)=\sum_{\substack{m \sim U \\ m\equiv0({\mathrm{mod}}~d)}} R(m);  f_d(\alpha)=\sum_{U<x\leq2U\atop x\equiv0\mod d}e(\alpha x^4);$$ $$g(\alpha)=\sum_{U<p\leq2U}e(\alpha p^4); h(\alpha)=\sum_{V<p\leq2V}e(\alpha p^4)$$
and then define
$$F(\alpha)=\sum_{d<D}\eta_df_d(\alpha),$$
$$v_{d}(n,\mathfrak{B})=\int_{\mathfrak{B}}f_{d}(\alpha)\overline{g(\alpha)}|g(\alpha)|^4|h(\alpha)|^{8}e(-\alpha n)d\alpha.$$
Naturally, $v_{d}(n,[0,1])=v_d(n)$, therefore
\begin{equation}
    \begin{split}
         |\sum_{d\leqslant D}\eta_dE_d(n)|&\leq \sum_{d\leqslant D}|\eta_d||v_d(n)-\frac{\mathfrak{S}_d(n)}{d} J(n)|\\&\leq \sum_{d\leqslant D}|\eta_d||v_d(n,\mathfrak{M})-\frac{\mathfrak{S}_d(n)}{d} J(n)|+\sum_{d\leqslant D}|\eta_d||v_d(n,\mathfrak{m})|.
    \end{split}
\end{equation}
Next, we calculate the second part on the right side of the inequality
\begin{equation}
    \begin{split}
        \sum_{d\leqslant D}|\eta_d||v_d(n,\mathfrak{m})|&=\int_{\mathfrak{m}}|F(\alpha)||g(\alpha)|^7|h(\alpha)|^{8}e(-\alpha n)d\alpha\\&\leq \bigg(\int_{\mathfrak{m}}|F(\alpha)|^2|g(\alpha)|^2|h(\alpha)|^{4}d\alpha\bigg)^{1/2}\bigg(\int_{\mathfrak{m}}|g(\alpha)|^{12}|h(\alpha)|^{12}d\alpha\bigg)^{1/2}.
    \end{split}
\end{equation}
Now we calculate the upper bound of $F(\alpha)$. Regarding the above $\mathfrak{M}'(q,a)$, $1\leq a\leq q\leq U^2$, we write $F(\alpha )$ as
$$F(\alpha)=\sum_{d\leq D}\eta_d\sum_{U/d<y\leq2U/d}e(\alpha d^4y^4).$$
By Dirichlet' s theorem on diophantine approximation, there are coprime integers $b$, $r$ with $r\leq16P^3d^{-3}$, $\left|d^{4}\alpha-\frac{b}{r}\right|\leq\frac{1}{16}r^{-1}d^{3}P^{-3}$. By Weyl' s inequality, when $r>U/d$
$$\sum_{U/d<y\leq2U/d}e(\alpha d^4y^4)\ll \bigg(\frac{U}{d}\bigg)^{7/8},$$
and when $r\leq U/d$
\begin{equation}
    \sum_{U/d<y\leq2U/d}e(\alpha d^4y^4)\ll r^{-\frac{1}{4}}\frac{U}{d}\bigg(1+\bigg(\frac{U}{d}\bigg)^4\bigg|\alpha d^4-\frac{b}{r}\bigg|\bigg)^{-\frac{1}{4}}+\bigg(\frac{P}{d}\bigg)^{\frac{1}{2}+\epsilon}.
\end{equation}
Furthermore, when
$$r\leq(U/d)^{7/8};\quad\left|\alpha d^{4}-\frac{b}{r}\right|\leq\frac{1}{r}\left(\frac{d}{U}\right)^{\frac{7}{2}},$$
it can also achieve a result of (3.6).
Hence
$$F(\alpha)\ll U^{7/8+\epsilon}D^{1/8}+U\sum_{d\in\mathcal{D}}d^{-1}r^{-\frac{1}{4}}\biggl(1+\biggl(\frac{U}{d}\biggr)^{4}\biggl|\alpha d^{4}-\frac{b}{r}\biggr|\biggr)^{-\frac{1}{4}},$$
where $\mathcal D $ represents the set of d that satisfies the condition. Compare the conditions of $q,a$ and $b,r$, we have
$$\bigg|\frac{b}{r}-\frac{ad^4}{q}\bigg|\leq \frac{1}{r}\bigg(\frac{d}{U}\bigg)^{\frac{7}{2}}+\frac{d^4}{qU^2}$$
i.e.
    $$\bigg|bq-ad^4r\bigg|\leq qd^{\frac{7}{2}}U^{-\frac{7}{2}}+rU^{-2}D^4\ll 1$$
where $D = N^\xi$, $0 < \xi < 9/25$, $U$ is large enough. Therefore $bq=ad^4r$, then $r=q/(q,d^4)$, by the trivial bound $(q,d^4)\leq (q,d)^4$
$$\sum_{d\in \mathcal D} d^{-1}r^{-\frac{1}{4}} \left( 1 + \left(\frac{U}{a}\right)^4 \left| \alpha d^4 - \frac{b}{r} \right| \right)^{-\frac{1}{4}} \leq q^{-\frac{1}{4}} \left( 1 + P^4 \left| \alpha - \frac{a}{q} \right| \right)^{-\frac{1}{4}} \sum_{d\leq D} \frac{(q,d)}{d}.$$
Thus
$$F(\alpha)\ll U^{\frac{7}{8}+\epsilon}D^{\frac{1}{8}}+q^{\epsilon-\frac{1}{4}}P(\log P)\left(1+P^{4}\bigg|\alpha-\frac{a}{q}\bigg|\right)^{-\frac{1}{4}}.$$
Let $\delta >0$ be so small that $D^{\frac{1}{8}}U^{\frac{7}{8}}\leq U^{\frac{23}{25}-2\delta}$, this is always possible because $D<U^{\frac{9}{25}}$. Let $\mathfrak{N}(q,a)$ denote the interval $|q\alpha-a|\le U^{5\delta-\frac{92}{25}}$, and $\mathfrak{N}$ be the union of all $\mathfrak{N}(q,a)$ with $1\le a\le q\le P^{\frac{8}{25}+5\delta}$, $(a,q)=1$. By the upper bound of $F(\alpha)$ obtained from the above, for $|F(\alpha)>U^{\frac{23}{25}}|$, then $\alpha\in \mathfrak{N}$. Defining $\Phi$ on $\mathfrak{N}$ as 
$$\Phi(\alpha)=q^{\epsilon-\frac{1}{4}}\Bigg(1+U^{4}\Bigg|\alpha-\frac{a}{q}\Bigg|\Bigg)^{-1/4}$$
Then we have
\begin{equation}
    \int_{\mathfrak{m}}|Fgh^2|^2d\alpha\ll U^{\frac{46}{25}-2\delta}\int_0^1|gh^2|^2d\alpha+U^2(\log U)^2\int_{\mathfrak{N}\cap\mathfrak{m}}|\Phi gh^2|^2d\alpha.
\end{equation}
By Hölder' s inequality
\begin{equation}
    \int_{\mathfrak{N}\cap\mathfrak{m}}|\Phi gh^2|^2d\alpha\leq\left(\int_{\mathfrak{N}\cap\mathfrak{m}}|\Phi|^{12}d\alpha\right)^{\frac{1}{6}}U_1^{\frac{1}{2}}U_2^{\frac{1}{3}}
\end{equation}
where
   $$ U_1 = \int_{0}^{1} |g(\alpha)|^2|h(\alpha)|^4 d\alpha; \quad U_2 = \int_{0}^{1} |g(\alpha)|^3|h(\alpha)|^6 d\alpha.$$
By simple calculations we obtain that
$$\int_{\mathfrak{N}\cap\mathfrak{m}}|\Phi(\alpha)|^{12}d\alpha\ll U^{-4}L^{-B}.$$
Then let $\mathfrak{m}_1(q,a)$ denote the interval $|q\alpha-a|\le U^{\frac{12}{5}}$, $\mathfrak{m}_1$ be the union of all $\mathfrak{m}_1(q,a)$ with $1\le a\le q\le U^{\frac{4}{5}}$, $(a,q)=1$, and $\mathfrak{m}_2(q,a)$ denote the interval $|q\alpha-a|\le V^{\frac{12}{5}}$, $\mathfrak{m}_2$ be the union of all $\mathfrak{m}_2(q,a)$ with $1\le a\le q\le V^{\frac{4}{5}}$, $(a,q)=1$. By Theorem 4.1, Lemma 6.3 of Vaughan\textsuperscript{\cite{7}}, and Theorem 2 of Vaughan\textsuperscript{\cite{11}}, we have
$$g(\alpha)\ll q^{-\frac14}U\bigg(1+U^4\bigg|\alpha-\frac aq\bigg|\bigg)^{-1}+U^{4/5+\epsilon}$$
where $\alpha\in \mathfrak{m}_1(q,a)$;
$$h(\alpha)\ll q^{-\frac14}V\bigg(1+V^4\bigg|\alpha-\frac aq\bigg|\bigg)^{-1}+V^{4/5+\epsilon}$$
where $\alpha\in \mathfrak{m}_2(q,a).$
Moreover, $|g(\alpha)|>U^{4/5+2\epsilon}$ implies $\alpha \in {\mathfrak{m}_1 }({\mathrm{mod}}~1)$, $|h(\alpha)|>V^{4/5+2\epsilon}$ implies $\alpha \in {\mathfrak{m}_2 }({\mathrm{mod}}~1)$.
\begin{equation}
    g(\alpha)\ll q^{-1/4}U+U^{1/2}\ll q^{-1/4}U
\end{equation}
where $\alpha \in {\mathfrak{m}_1 }({\mathrm{mod}}~1)$;
\begin{equation}
    h(\alpha)\ll q^{-1/4}V+V^{1/2}\ll q^{-1/4}V
\end{equation}
where $\alpha \in {\mathfrak{m}_2 }({\mathrm{mod}}~1).$
Then define $\Psi(\alpha)$ on $\mathfrak{m}_1$, $\Psi^{*}(\alpha)$ on $\mathfrak{m}_2$ as
\begin{equation}
    \Psi(\alpha)=q^{-\frac{1}{4}}\left(1+U^4\bigg|\alpha-\frac{a}{q}\bigg|\right)^{-1};\Psi^{*}(\alpha)=q^{-\frac{1}{4}}\left(1+V^4\bigg|\alpha-\frac{a}{q}\bigg|\right)^{-1},
\end{equation}
Subsequently
\begin{align*}
   \int_{0}^{1} |g(\alpha)|^3|h(\alpha)|^6 d\alpha\ll &\: U^{\frac45}\int_{0}^{1} |g(\alpha)|^2|h(\alpha)|^6 d\alpha+U\int_{\mathfrak{m}_1}|\Psi(\alpha)|^2|g(\alpha)|^2|h(\alpha)|^6 d\alpha\\\ll &\: U^{\frac45} V^{\frac85}\int_{0}^{1} |g(\alpha)|^2|h(\alpha)|^4d\alpha+ U^{\frac45}V^2\int_{\mathfrak{m}_2}|\Psi^*(\alpha)|^2|g(\alpha)|^2|h(\alpha)|^4 d\alpha\\&+UV^{\frac85}\int_{\mathfrak{m}_1}|\Psi(\alpha)|^2|g(\alpha)|^2|h(\alpha)|^4 d\alpha\\&+UV^2\int_{\mathfrak{m}_1\cap\mathfrak{m}_2}|\Psi(\alpha)|^2|\Psi^*(\alpha)|^2|g(\alpha)|^2|h(\alpha)|^4 d\alpha.
\end{align*}
By (3.9), (3.10), (3.11) we can obtain that
\begin{equation}
    U_2=\int_{0}^{1} |g(\alpha)|^3|h(\alpha)|^6 d\alpha\ll U^{5+13/40}.
\end{equation}
Since the upper bound of $U_1$ is known, we can combine(3.7), (3.8), (3.12) to obtain that
$$ \int_{\mathfrak{m}}|Fgh^2|^2d\alpha\ll  U^{4+193/200}L^{-B/6}.$$
Similar to the method mentioned above, we can obtain
$$\int_{\mathfrak{m}}|g(\alpha)|^{12}|h(\alpha)|^{12}d\alpha\ll U^{16+29/40}$$
Naturally, by (3.5)
$$ \sum_{d\leqslant D}|\eta_d||v_d(n,\mathfrak{m})|\ll U^{11}L^{-B/12}$$
Now we calculate the first term on the right side of the inequality (3.4). For $\alpha=a/q+\beta\in \mathfrak{M}$, we define 
$$f_{d}^{*}(\alpha)=\frac{S(q,ad^{4})}{qd}\Phi(\lambda,U),\quad g^{*}(\alpha)=\frac{C(q,a)}{\varphi(q)}\Psi(\lambda,U),\quad h^{*}(\alpha)=\frac{C(q,a)}{\varphi(q)}\Psi(\lambda,V).$$
By Lemma 7.15 of Hua\textsuperscript{\cite{12}}
$$g(\alpha)-g^*(\alpha)\ll U\exp(-c_1\sqrt{L}),\quad h(\alpha)-h^*(\alpha)\ll V\exp(-c_1\sqrt{L})$$
where $c_1>0$ is  a absolute constant.By Theorem 4.1 of Vaughan\textsuperscript{\cite{11}}
$$f_{d}(\alpha)-f_{d}^{*}(\alpha)\ll q^{1/2+\epsilon}(1+U^{4}|\lambda|)\ll L^{2B}.$$
Then for $\alpha\in \mathfrak M$, by trivial estimate $\left|S(q,ad^{3})\right|\leq q$, we have
$$f_{d}\overline{g}|g|^{6}|h|^{8}-f_{d}^{*}\overline{g^{*}}|g^{*}|^{6}|h^{*}|^{8}\ll d^{-1}U^{8}V^{8}L^{-13B},$$
and 
\begin{equation}
    \sum_{d\leq D}|\eta_{d}||\nu_{d}(n,\mathfrak{M})-\nu_{d}^{*}(n)|\ll U^4V^{8}L^{-A},
\end{equation}
where
$$\nu_{d}^{*}(n)=\int_{\mathfrak{M}}f_{d}^{*}\overline{{g^{*}}}|g^{*}|^{6}|h^{*}|^{8}e(-n\alpha)d\alpha,$$
and we can find that
\begin{equation}
    \nu_{d}^{*}(n)=\frac{1}{d}J(n)\sum_{q\leq L^{B}}T_{d}(n,q).
\end{equation}
Then by $|C(q,a)|\ll q^{3/4+\epsilon}$
$$|T_{d}(n,q)|\leq\sum_{a=1}^{q}\frac{|S(q,ad^{4})||C(q,a)|^{15}}{q\varphi^{15}(q)}\ll\sum_{a=1}^{q}\frac{q^{49/4+\epsilon}}{q\varphi^{15}(q)}\ll q^{-2}.$$
Hence the series
$$\mathfrak{S}_{d}(n)=\sum_{q=1}^{\infty} T_{d}(n,q)$$
converges absolutely, this proves (i), and
$$\sum_{q\leq L^{B}} T_{d}(n,q)=\mathfrak{S}_{d}(n)+O(L^{-B}).$$
By (3.14)
$$\nu_d^*(n)=\frac{\mathfrak{S}_d(n)}{d}J(n)+O\left(\frac{U^4V^8}{dL^B}\right),$$
combining with (3.13), we can obtain
$$\sum_{d\leqslant D}|\eta_d||v_d(n,\mathfrak{M})-\frac{\mathfrak{S}_d(n)}{d} J(n)|\ll U^4V^8L^{-A}.$$
This proves (iii), which also proves Lemma 2.\\

\textbf{Lemma 3.}\textit{ For $(d,6)=1$, we get }
\begin{align*}
    \begin{split}
        \mathfrak{S}_d(n)=&\{1+T_1(n,2)+T_1(n,2^2)+T_1(n,2^3)+T_1(n,2^4)\}\\&\times\prod_{\begin{array}{c}p\nmid d\\p\neq2\end{array}}\{1+T_1(n,p)\}\prod_{p\mid d}\{1+T_p(n,p)\}.
    \end{split}
\end{align*}
\textbf{Proof.} According to the definition in (2.3), it can be obtained through a simple proof $\mathfrak{S}_d (n)$ which is an integral function, then there is
\begin{align*}
    \begin{split}
        \mathfrak{S}_d(n)&=\prod_{p}\{1+T_p(n,p)+T_p(n,p^2)+\cdots \}\\&=\prod_{p\nmid d}\{1+T_1(n,p)+T_1(n,p^2)+\cdots\}\\&\times\prod_{p\mid d}\{1+T_p(n,p)+T_p(n,p^2)+\cdots\}.
    \end{split}
\end{align*}
For $p\nmid d$, we have
\begin{align*}
    \begin{split}
        T_d(n,p)&=\sum_{\overset{a=1}{(a,p)=1}}^p\frac{S(p,ad^3)C^3(p,a)\overline{C^4(p,a)}}{p\varphi^7(p)}e\left(-\frac{an}p\right)\\&=\sum_{\overset{a=1}{(a,p)=1}}^p\frac{S(p,a)C^3(p,a)\overline{C^4(p,a)}}{p\varphi^7(p)}e\left(-\frac{an}p\right)=T_1(n,p).
    \end{split}
\end{align*}
According to Lemma 4 in Hua\textsuperscript{\cite{6}}, $C(p^t,a)=0$($p\neq2$ and $t\geq2$ or $p=2,t\geq5$), substituting back to the original equation yields the Lemma 3.\\

\noindent{\textbf{Lemma 4}.} \textit{Define $K(n,p)$ as the number of solutions to the following equation}

$$y_2^4 + \cdots + y_8^4 - y_9^4 - \cdots - y_{16}^4 \equiv n (\text{mod}~p)$$
\textit{where $1\leq y_i\leq p-1$ $(2\leq i\leq 16)$, then we get}
$$pK(n,p) = (p - 1)^{15} + E,$$
\textit{where }

$$E=\begin{cases}\:-(p-1),&p\mid n,\\\:1,&p\nmid n,\end{cases}~~~~~~~~~~~~~~\textit{if}~p\equiv3({\mathrm{mod}}~4);$$
\textit{and}
$$|E|\leqslant(3\sqrt{p}+1)^{{13}}(p-1)(3p+1),~~~\textit{if}~p\equiv1({\mathrm{mod}}~4).$$
\textit{It follows that $K(n,p)>0$ for} $p\geq17$.\\
\textbf{Proof.} From the definition of $K(n,p)$, we have
$$pK(n,p)=\sum_{a=1}^{p}C^{7}(p,a)\overline{C^{8}(p,a)}e\left(-\frac{an}{p}\right)=(p-1)^{15}+E,$$
where
$$E=\sum_{a=1}^{p-1}C^{7}(p,a)\overline{C^{8}(p,a)}e\biggl(-\frac{an}{p}\biggr).$$
When $p\equiv3(\text{mod }4)$ and $(p,a)=1$, by Lemma 4.3 in Vaughan\textsuperscript{\cite{7}}, we get $S(p,a)=0$ thus $C(p,a)=-1$, therefore
$$E=-\sum_{a=1}^{p-1}e\left(-\frac{an}{p}\right)=\begin{cases}-(p-1),&p\mid n,\\1,&p\nmid n.\end{cases}$$
and when $p\equiv1(\text{mod }4)$, by Lemma 4.3 of Vaughan\textsuperscript{\cite{7}} again, we have $|C(p,a)|\leqslant3\sqrt{p}+1$, and for

$$\begin{aligned}\sum_{a=1}^{p-1}|C(p,a)|^2=\sum_{a=1}^p|C(p,a)|^2-(p-1)^2.\end{aligned}$$

Obviously,  $\sum_{a=1}^p|C(p,a)|^2$ can be expressed as  $p$ times the number of solutions to equation $x^4\equiv y^{4}({\mathrm{mod~}}p)$, $1\leqslant x$, $y\leqslant p- 1$.
For $p\equiv1({\mathrm{mod~}}4)$, $$\sum_{a=1}^p|C(p,a)|^2=4p(p-1),$$
then
$$\sum_{a=1}^{p-1}|C(p,a)|^2=(3p+1)(p-1),$$
thus
$$\sum_{a=1}^{p-1}|C^7(p,a)\overline{C^8(p,a)}|\leqslant(3\sqrt{p}+1)^{13}(p-1)(3p+1).$$
 A simple calculation shows that $K(n,p)>0$ for $p\geq 15$, $p\equiv1(\text{mod }4)$and for all $p\equiv3(\text{mod }4)$, therefore we have chosen $p\geq17$. This proves Lemma 4.\\

From the  similar method of Lemma 4, we can get Lemma 5.\\

\noindent{ \textbf{Lemma 5.}} \textit{For $i = 1,2$, let $H(n, p^i) $ denote the number of solution of}

$$x^4 + y_2^4 + \cdots + y_8^4 - y_9^4 - \cdots - y_{16}^4 \equiv n \pmod{p^i}$$
\textit{where $1\leq x \leq p^i, 1 \leq y_j < p^i $ and $ (y_j, p) = 1$. Thus}
$$pH(n, p) = p(p-1)^{15} + E^*,$$
\textit{where}
$$ E^* = 0, ~~~~~\textit{ if}~p\equiv3(\text{mod }4);$$
and
$$|E^*| \leq 3\sqrt{p}(3\sqrt{p}+1)^{13}(p-1)(3p+1),~~\textit{if}~p\equiv1(\text{mod }4).$$
\textit{It also follows that $H(n, p)>0 $ for $p>17$.}\\
\\
\noindent\textbf{Lemma 6.} \textit{Definne the functions }
$$\phi_{0}(u)=\frac{2e^{\gamma_{0}}}{u}\mathrm{log}(u-1)\quad \text{and}\quad\phi_{1}(u)=\frac{2e^{\gamma_{0}}}{u},$$
\textit{where $2\leq u\leq 3$. Suppose $\omega(d)$ is a multiplicative function of $d$ satisfying the conditions }
$$0\leqslant\omega(p)<p\quad \text{and}\quad\omega(p^l)=1+O(p^{-1}),$$
\textit{for each prime number $p$ and natural number $l$. Let $X$ be a real number with $X>3$, for $r (x)$ be a non-negative arithmetical function, we define }
$$E_{d}=\sum_{\begin{smallmatrix}P<x<2P\\x\equiv0\pmod{d}\end{smallmatrix}}r(x)-\frac{\omega(d)}{d}X.$$
\textit{Let $U,V $ and $ z$ be positive real parameters satisfying the inequality }
$$2\leqslant\frac{\log(UV)}{\log z}\leqslant3.$$
\textit{For any sequences $\{a_m\}$ and $\{b_k\}$ with}
$$|a_m|\leqslant1\quad and\quad|b_k|\leqslant1,$$
\textit{one has}
$$\sum_{1\leq m\leq U}a_{m}\sum_{1\leq k\leq V}b_{k}E_{mk}\ll X(\log X)^{-2}.$$
\textit{Then, we write}
$$W(z)=\prod_{p<z}\:(1-\omega(p)/p),$$
\textit{one has the lower bound}
$$
\sum_{\substack{P \leq x < 2P \\ (x, \Pi(z)) = 1}} r(x) > X W(z) \left( \phi_0 \left( \frac{\log(UV)}{\log z} \right) + O \left( (\log \log X)^{-1/50} \right) \right),$$
\textit{and also the upper bound}
$$\sum_{\substack{P \leq x < 2P \\ (x, \Pi(z)) = 1}} r(x) < X W(z) \left( \phi_1 \left( \frac{\log(UV)}{\log z} \right) + O \left( (\log \log X)^{-1/50} \right) \right).$$
\textbf{Proof.} See [8], Lemma 9.1.
\\

\noindent\textbf{Lemma 7.} \textit{(Mertens' theorem)For prime number $p$, $x\geq e$, positive integers $k$, $l$ with $(k,l)=1$ and $k\leq \ln^A{x}$ for any positive integer A, then}
$$\prod_{\substack{p\leq x\\p\equiv l({\mathrm{mod~}k})}}\left(1-\frac1p\right)=\frac{e^{-\frac1{\varphi(k)}\left[\gamma+\ln\frac{\varphi(k)}k\right]}}{\ln^{\frac1{\varphi(k)}}x}\left\{1+O\left(e^{-c\ln^{\frac35}x}\right)\right\}$$
\textit{where $c$ is a positive absolute constant, $\gamma$ is the Euler's constant, and the constant in $O$ is independent of $x$.}
\textbf{Proof.} See [15], corollary of Theorem 429.
\\

\noindent\textbf{Lemma 8.} \textit{For $N/9 < \ell \leqslant N$, define $r(\ell)$  as the number of $\ell$ can be expressed in the form of (1.2) with
$$p_{1}, p_{2},p_{3}, p_{4}\sim U, \qquad p_{5}, p_{6}, p_{7}, p_{8}\sim V \, .$$
Then we have
$$\sum_{N/9<\ell\leqslant N}r^2(\ell)<bU^2V^8L^{-14},$$
where $b= 80947432211.141$.}\\
\textbf{Proof.} Denote $\mathcal{B}$  as the set of all prime numbers greater than 17, and define
\begin{equation*}
    P(z) = \prod_{\substack{p < z \\ p \in \mathcal{B}}} p.
    \label{eq:placeholder_label}
\end{equation*}
To prove Lemma 8, we hope to obtain an appropriate upper bound for the following equation
$$\sum_{\substack{m\sim U \\ \atop(m,P(z))=1}}R(m).$$
For this, we apply Lemma 6, there is no prime divisor beyond $\mathcal{B}$ for $d$, then we define
$$\omega(d)=\mathfrak{S}_d(n)/\mathfrak{S}_1(n).$$
Especially by Lemma 3, for $p\in \mathcal{B}$ we have
\begin{equation*}
    \omega(p)=\frac{1+T_{p}(n,p)}{1+T_{1}(n,p)}.
\end{equation*}
By (2.3), we get
$$\begin{aligned}1+T_p(n,p)=\sum_{a=1}^p\frac{\overline{C(p,a)}|C(p,a)|^{14}}{(p-1)^{15}}e\Bigg(-\frac{an}{p}\Bigg)=p\frac{K(n,p)}{(p-1)^{15}}.\end{aligned}$$
hence
$$\omega(p)=\frac{K(n,p)}{H(n,p)}.$$
By Lemma 4 and Lemma 5, through simple calculations, for all $p\in \mathcal{B}$ and the positive integer $l$, we have
$$0\leqslant\omega(p)<p,\quad\omega(p^l)=1+O(p^{-1}).$$

Now let $X=\mathfrak{S}_1(n)J(n)$, then by Lemma 2, we have
$$\sum_{\substack{m\sim U\\m\equiv 0(\mod d)}}R(m)=\frac{\omega(d)}{d}X+E_d(n).$$
Suppose $U_{0}\geqslant 1,\ V_{0}\geqslant 1,\ U_{0}V_{0}=D=N^{9(1-\epsilon)/100}$ and $z=D^{1/2}$, then
$$\frac{\log(U_{0}V_{0})}{\log z}=2.$$
For any sequence of numbers $\{a_m\}$,  $\{b_k\}$ with
$$|a_m|\leq 1~ ~, ~~|b_k|\leq1,$$
by Lemma 2, we have
$$\sum_{1\leqslant m\leqslant U}a_m\sum_{1\leqslant k\leqslant V}b_kE_{mk}\ll\sum_{1\leqslant d\leqslant D}\tau(d)E_d\ll U^4V^8L^{-A}.$$
By the upper bound in Lemma 6, we get
$$\sum\limits_{(m,P(z))=1}R(m)<e^{\gamma}(1+\varepsilon)J(n)\mathfrak{S}_{1}(n)W(z),$$
where $\gamma$ denotes Euler's constant, and
$$W(z)=\prod_{\stackrel{p\in\mathcal{B}}{p<z}}\biggl(1-\frac{\omega(p)}{p}\biggr).$$
Actually, we have estimated the upper bound of $J (n)$ in Lemma 2. Next, we estimate the remaining part,
\begin{equation}
    \begin{split}
&\mathfrak{S}_{1}(n)W(z)\\
=&\{1+T_1(n,2)+T_1(n,2^2)+T_1(n,2^3)+T_1(n,2^4)\}\\
&\times(1+T_1(n,3))(1+T_1(n,5))(1+T_1(n,7))(1+T_1(n,11))(1+T_1(n,13))\\
&\times\prod_{17\leqslant p<N^{9(1-\varepsilon)/200}}(1+T_1(n,p))\biggl(1-\frac{K(n,p)}{H(n,p)}\biggr)\\
&\times\prod_{p\geqslant N^{9(1-\varepsilon)/200}}(1+T_1(n,p))\\=&\{1+T_1(n,2)+T_1(n,2^2)+T_1(n,2^3)+T_1(n,2^4)\}\\
       &\times(1+T_1(n,3))(1+T_1(n,5))(1+T_1(n,7))(1+T_1(n,11))(1+T_1(n,13))\\
       &\times\prod_{17\leqslant p<N^{(1-\varepsilon)/200}}\left(1-\frac{1}{p}\right)\left(1+\frac{E^{*}-E}{\left(p-1\right)^{16}}\right)\\
       &\times\prod_{p\geqslant N^{(1-\varepsilon)/200}}\left(1+\frac{E^{*}}{p(p-1)^{15}}\right).
    \end{split}
\end{equation}

For $n$ in the definition of $R(m)$, define $\rho(n)$ as the number of solutions to the following equation
$$n=p_1^4+\cdots+p_8^4-p_9^4-\cdots-p_{16}^4,$$
where
$$p_1,p_2,p_3,p_4,p_9,p_{10},p_{11},p_{12}\sim U,\quad p_5,p_6,p_7,p_{8},p_{13},p_{14},p_{15},p_{16}\sim V.$$
We can naturally obtain
$$\sum_{N/9<\ell\leqslant N}r^2(\ell)\leqslant\rho(0)$$
and
$$\rho(n)\leq \sum_{m\sim U\atop(m,P(z))=1}R(m).$$

Now we are calculating the various parts in (3.15). Since we are actually just trying to obtain a suitable upper bound for $\rho(0)$, thus we assume $n=0$ in the following text.\\

According to the definition (2.3), (2.4), by simple calculation we have $S(2,1)=0$, then
$$T_{1}(n,2)=\sum_{\overset{a=1}{(a,2)=1}}^2\frac{S(2,a)C^7(2,a)\overline{C^8(2,a)}}{2\varphi^{16}(2)}e\biggl(-\frac{an}{2}\biggr)=0.$$

For $q=4$, we have
$$T_1(n,4)=\sum_{\overset{a=1}{(a,4)=1}}^4\frac{S(4,a)C^7(4,a)\overline{C^8(4,a)}}{4\varphi^{16}(4)}e\biggl(-\frac{an}{4}\biggr),$$
where $$S(4,a)=\sum_{m=1}^4e\biggl(\frac{a m^4}{4}\biggr)=2+2e\biggl(\frac{a}{4}\biggl), C(4,a)=\sum_{(m,4)=1}^4e\biggl(\frac{a m^4}{4}\biggr)=2e\biggl(\frac{a}{4}\biggl).$$
Thus
\begin{align*}
    \begin{split}
        T_1(n,4)&=\sum_{\overset{a=1}{(a,4)=1}}^4\frac{(2+2e(\frac{a}{4}))(2e(\frac{a}{4}))^7\overline{2e(\frac{a}{4})}^8}{4\varphi^{13}(4)}e\biggl(-\frac{a\times0}{4}\biggr)\\&=\frac{S(4,1)C^{7}(4,1)+S(4,3)C^{7}(4,3)}{2^{15}}=1.
    \end{split}
\end{align*}
Similarly, we have
$$T_1(n,8)=\sum_{\overset{a=1}{(a,8)=1}}^8\frac{1}{2}\bigg(1+e\bigg(\frac{7a}{8}\bigg)\bigg)=2$$
and $$T_1(n,16)=\sum_{\overset{a=1}{(a,16)=1}}^{16}\frac{1}{2}\bigg(1+e\bigg(\frac{15a}{16}\bigg)\bigg)=4.$$
Thus
$$1+T_1(n,2)+T_1(n,2^2)+T_1(n,2^3)+T_1(n,2^4)=8.$$

Also, for $q=3$, we have $$T_1(n,3)=\sum_{\overset{a=1}{(a,3)=1}}^3\frac{S(3,a)C^7(3,a)\overline{C^8(3,a)}}{3\varphi^{15}(3)}e\biggl(-\frac{an}{3}\biggr),$$
where$$S(3,a)=\sum_{m=1}^3e\biggl(\frac{a m^4}{3}\biggr)=1+2e\bigg(\frac{a}{3}\bigg), C(4,a)=\sum_{(m,4)=1}^4e\biggl(\frac{a m^4}{4}\biggr)=2e\bigg(\frac{a}{3}\bigg).$$
then
\begin{align*}
    \begin{split}
        T_1(n,3)&=\sum_{\overset{a=1}{(a,3)=1}}^3\frac{(1+2e(\frac{a}{3}))(2e(\frac{a}{3}))^7\overline{(2e(\frac{a}{3}))^8}}{3\varphi^{15}(3)}e\biggl(-\frac{a\times 0}{3}\biggr)=1.
    \end{split}
\end{align*}
Likewise, by a series of calculations we have $T_1(n,5)=3$, $T_1(n,7)=\frac{128}{3^{14}}$, $T_1(n,11)=\frac{2187}{6103515625 }$ and $T_1(n,13)\leq 0.015$ .

Furthermore, for
$$\prod_{17\leqslant p<N^{9(1-\varepsilon)/200}}\left(1+\frac{E^{*}-E}{\left(p-1\right)^{16}}\right),$$
by Lemma 4 and Lemma 5, for $p\equiv3(\text{mod }4)$ we have
$$|E-E^*|\leq p-1,$$
then
\begin{align*}
    \begin{split}
        \prod_{\overset{17\leqslant p<N^{9(1-\varepsilon)/200}}{p\equiv3(\text{mod}4)}}\left(1+\frac{E^{*}-E}{\left(p-1\right)^{16}}\right)&\leq \prod_{\overset{p\geqslant17}{p\equiv3(\mathrm{mod~4})}}\left(1+\frac{1}{\left(p-1\right)^{15}}\right)\\&<1\:+\:\sum_{n\geqslant19}\:\frac{1}{n^{15}}<1\:+\:19^{-14},
    \end{split}
\end{align*}
and for $p\equiv1(\text{mod }4)$, we have
$$|E-E^*|\leq (3\sqrt{p}+1)^{14}(p-1)(3p+1).$$
Now define $\delta_p$ as
$$\left|\frac{E^*-E}{(p-1)^{16}}\right|\leq\frac{(3\sqrt{p}+1)^{14}(3p+1)}{(p-1)^{15}}=\delta_p.$$
By numerical calculation
$$\sum_{p \geq 17}\delta_p<0.03,$$
then by $1+x\leq e^x$
$$\prod_{\overset{17\leq p \leq N^{{9(1-\epsilon)}/{200}}}{p\equiv1(\text{mod}4)}}\left(1+\frac{E^{*}-E}{\left(p-1\right)^{16}}\right)\leq \prod_{\overset{p \geq 17}{p\equiv1(\text{mod}4)}}e^{\frac{E^{*}-E}{\left(p-1\right)^{16}}}\leq e^{\sum_{p\geq 17} {\delta_p}}\leq 1.1.
$$
Therefore we have
\begin{equation*}
    \prod_{13\leqslant p<N^{(9(1-\varepsilon)/200}}\left(1+\frac{E^{*}-E}{\left(p-1\right)^{16}}\right)<1.1
\end{equation*}

Furthermore, by Lemma 7, we have the estimate
$$\prod_{17\leqslant p<N^{(1-\varepsilon)/32}}\left(1-\frac{1}{p}\right)<231.713 e^{-\gamma}(1+\varepsilon)L^{-1}.$$

It is not difficult to find that for $p$ large enough. 
$$\frac{|E^*|}{p(p-1)^{15}}<\frac1{p^2}$$
 Hence
$$\prod_{p\geqslant N^{9(1-\varepsilon)/200}}\left(1+\frac{E^*}{p(p-1)^{15}}\right)\leqslant\prod_{p\geqslant N^{9(1-\varepsilon)/200}}\left(1+\frac{1}{p^2}\right)<1+\varepsilon.$$\\

To sum up,
$$\mathfrak{S}_{1}(n)W(z)<16557.733e^{-\gamma}L^{-1}$$
then
$$\sum_{m\sim U\atop(m,P(z))=1}R(m)<bU^4V^8L^{-16},$$
where $b=80947432211.141$.

Naturally, for $n=0$,
$$\rho(0)\leq \sum_{m\sim U\atop(m,P(z))=1}R(m)$$
and
$$\sum_{N/9<\ell\leqslant N}r^2(\ell)\leq \rho(0),$$
therefore
$$\sum_{N/9<\ell\leqslant N}r^2(\ell)<bU^4V^8L^{-16}.$$
We have proved Lemma 8.
\section{Proof of the Theorem}
In this part, we will prove Theorem.\\

By the prime number theorem, for
$$U=\left(\frac{N}{64(1+\delta_0)}\right)^{1/4},\quad V=U^{7/8},$$
we have
\begin{align}
    \begin{split}
        \sum_{N/9<\ell\leqslant N}r(\ell)&\geqslant\sum_{p_1\sim U}1\sum_{p_2\sim U}1\sum_{p_3\sim U}1\sum_{p_4\sim U}1\sum_{p_5\sim V}1\sum_{p_6\sim V}1\sum_{p_7\sim V}\sum_{p_8\sim V}\\&\geqslant(1-\varepsilon)\frac{U^4V^4}{\log^4U\log^4V}\geqslant\left(\frac{128}{7}\right)^4(1-\varepsilon)U^4V^4L^{-8}.
    \end{split}
\end{align}
then by Cauchy's inequality and Lemma 8, we have
\begin{align}
    \begin{split}
        \left\{\sum_{N/9<\ell\leqslant N}r(\ell)\right\}^{2}&\leqslant\left\{\sum_{\substack{N/9<\ell\leqslant N\\r(\ell)>0}}1\right\}\left\{\sum_{N/9<\ell\leqslant N}r^{2}(\ell)\right\}\\&\leqslant bU^4V^8L^{-16}\left\{\sum_{N/9<\ell\leqslant N}1\right\}.
    \end{split}
\end{align}
From (4.1), (4.2), we have
$$\sum_{N/9<\ell\leqslant N\atop r(\ell)>0}1>\frac{(1-\varepsilon)^{2}}{b}\biggl(\frac{128}{7}\biggr)^{8}U^{4}>\frac{1}{414.465}N.$$
The proof of our Theorem is now complete.

\hfill
$\square$


\vskip 3mm




\vskip 8mm

\begin{thebibliography}{150}

\bibitem{1}
Brüdern, A sieve approach to the Waring–Goldbach problem (I): sums of four cubes, Ann. Sci. Éc. Norm. Super. 28 (1995), 461-476.


\bibitem{2} X.M. Ren, Density of integers that are the sum of four cubes of primes, Chin. Ann. Math. Ser. B 22 (2001), 233-242.
\bibitem{3}X.M. Ren, Sums of four cubes of primes, J. Number Theory 98 (2003), 156-171.
\bibitem{4}Z.X. Liu, Density of the sums of four cubes of primes, J. Number Theory 132 (2012), 735-747.
\bibitem{5}K.F. Roth, On Waring’s problem for cubes, Proc. London Math. Soc. 53 (2) (1951), 268-279.
\bibitem{6} L.K. Hua, Some results in additive prime number theory, Quart. J. Math. (Oxford) 9 (1938), 68-80.
\bibitem{7} R.C. Vaughan, The Hardy–Littlewood Method, Cambridge University Press, Cambridge, 1981.
\bibitem{8} K. Kawada, T.D. Wooley, On the Waring–Goldbach problem for forth and fifth powers, Proc. London Math. Soc. 83 (3) (2001), 1-50.
\bibitem{9}J.Y. Liu, The Goldbach-Vinogradov theorem with three primes in a thin subset, Chinese Annals Math. 19(1998), 479-488.
\bibitem{10}K. Kawada, Note on the sum of cubes of primes and an almost prime, Arch. Math. 69 (1997), 13-19.
 \bibitem{11}R.C. VAUGHAN, Some remarks on Weyl sums (Topics in classical number theory, Colloq. Math. Soc. J. Bolyai, North Holland, Amsterdam. 34, 1984, 1585-1602). 
\bibitem{12}L.K. Hua, Additive theory of prime numbers (in Chinese) [M], Science Press, Beijing, 1957; English version, Amer. Math. Soc. Rhode Island, 1965.
\bibitem{13}W.J. Ellison, Waring’s problem. Amer. Math. Monthly, (1971), 78:10-36. 
\bibitem{14} Tom M. Apostol, Introduction to Analytic Number Theory [M].  New York: Springer-Verlag, 1976.
\bibitem{15}G.H. Hardy, and E.M. Wright, An Introduction to the Theory of Numbers [M]. Oxford University Press, 6th edition, 2008.
\end{thebibliography}
\end{document}